\documentstyle{article}
\date{}
\title{Hyperbolic Orbits for a Class of Singular
 Hamiltonian Systems\footnote{Supported partially by NSF (11071175) of China.}}
\author{{Donglun Wu and Shiqing Zhang}\\
{\small Yangtze Center of Mathematics and College of Mathematics, Sichuan University,}\\
{\small Chengdu 610064, People's Republic of China}}
\begin{document}
\maketitle \large \baselineskip 14pt

\begin{quote}

{\bf Abstract}\ \ The existence of hyperbolic orbits is proved for a
class of singular Hamiltonian systems $\ddot{u}(t)+\nabla V(u(t))=0$
by taking limit for a sequence of periodic solutions which are the
variational minimizers of Lagrangian actions.



\end{quote}

\section{Introduction and Main Results}

\ \ \ \ \ \ In this paper, we consider the following second order
Hamiltonian systems
\begin{equation}
   \ddot{u}(t)+\nabla V(u(t))=0\label{1}
\end{equation}
with
\begin{equation}
   \frac{1}{2}|\dot{u}(t)|^{2}+ V(u(t))=H\label{2}.
\end{equation}
where $u\in(R^{1},R^{N})$, $V\in C^{1}(R^{N}\setminus\{0\},R^{1})$
has a singularity at the origin. Subsequently, $\nabla V(x)$ denotes
the gradient with respect to the $x$ variable,
$(\cdot,\cdot):R^{N}\times R^{N} \rightarrow  R$\ denotes the
standard Euclidean inner product in $R^{N}$ and $\mid \cdot\mid$ is
the induced norm.

The periodic and homoclinic orbits of Hamiltonian systems have
been studied by many mathematicians [1-12, 14-17] and the
references therein. Specially, the n-body-type problem is a
Hamiltonian system which has attracted some mathematicians
 to use variational methods to study the parabolic
and hyperbolic orbits \cite{2,3,7,16}. Referring to the two-body
problem, with a center potential $V(x)=-\displaystyle\frac{1}{|x|}$,
it is well known that
\begin{eqnarray*}
&&(\mbox{i}).\ \mbox{If H$<$0, the solution for systems (\ref{1}) and (\ref{2}) is an elliptic orbit;}\nonumber\\
&&(\mbox{ii}).\ \mbox{If H=0, the solution for systems (\ref{1})
and (\ref{2}) is a parabolic orbit;}\nonumber\\
&&(\mbox{iii}).\ \mbox{If H$>$0, the solution for systems
(\ref{1}) and (\ref{2}) is a hyperbolic orbit.}
\end{eqnarray*}
Where the parabolic and hyperbolic orbits mean:

\vspace{0.3cm}\textbf{Definition 1.0.([3])}\ If when $|t|\rightarrow
+\infty$ we have
$$|u(t)|\rightarrow +\infty,|\dot{u}(t)|\rightarrow 0,$$
then we call $u(t)$ is a parabolic orbit;

 If when $|t|\rightarrow +\infty$, we have

$$|u(t)|\rightarrow +\infty,|\dot{u}(t)|>0,$$

then we call $u(t)$ is a hyperbolic orbit.

 In the above two cases, the parabolic and hyperbolic
orbits are all called hyperbolic-like orbits by Felmer and Tanaka
in \cite{3}. Subsequently, an orbit is said to be a parabolic or a
hyperbolic orbit, if it satisfies
\begin{eqnarray}
|u(t)|\rightarrow\infty\ \ \ \mbox{as}\ \ \
t\rightarrow\pm\infty.\label{5}
\end{eqnarray}
In 2000, for $N=2$, Felmer and Tanaka proved that

\vspace{0.3cm} \textbf{Theorem 1.1(See\cite{3}).}\ {\em Assume that $N=2$ and the following conditions hold\\

$(A_{1})$ $V\in C^{1}(R^{N}\setminus\{0\},R^{1})$,

$(A_{2})$ $V(x)\leq0$ for all $x\in R^{N}\setminus\{0\}$,

$(A_{3})$ there are constants $\zeta>2$, $\rho>0$ and $d_{0}>0$ such
that
\begin{eqnarray*}
&&(i).\ -V(x)\geq \frac{d_{0}}{|x|^{\zeta}}\ \ \ \mbox{for}\ \ 0<|x|\leq \rho,\nonumber\\
&&(ii).\ (x,\nabla V(x))+2V(x)\rightarrow +\infty\ \mbox{as}\
|x|\rightarrow0.\nonumber
\end{eqnarray*}

$(A_{4})$  there exist $\beta>2$ and $C_{0}>0$ such that
\begin{eqnarray*}
-V(x)\leq\frac{C_{0}}{|x|^{\beta}}\ \ \ \mbox{and}\ \ \ |\nabla
V(x)|\leq\frac{C_{0}}{|x|^{\beta+1}}\ \ \ \mbox{for}\ \ |x|\geq1.
\end{eqnarray*}

Then for any given $H>0$, $\theta_{+}$, $\theta_{-}\in R$ with
$\theta_{+}-\theta_{-}>\pi$, there exists a solution
$u(t)=r(t)(\cos\theta(t),\sin\theta(t))$ of $(1)-(2)$such that
$\theta(t)\rightarrow\theta_{\pm}$ as
$t\rightarrow\pm\infty$.}\\

For $N\geq3$, they proved that

\vspace{0.3cm}{\bf Theorem 1.2(See\cite{3})}\ \ {\em Assume
$N\geq3$ and $(A_{1})-(A_{4})$ hold. Then for any given $H>0$ and
$\theta_{+}\neq-\theta_{-}$, there exists a solution $u(t)$ of
$(1)-(2)$ such that}
\begin{eqnarray*}
\lim_{t\rightarrow\pm\infty}\frac{u(t)}{|u(t)|}=\theta_{\pm},
\end{eqnarray*}
where $\theta_{+}$, $\theta_{-}\in S^{N-1}=\{x\in R^{N}||x|=1\}$ are
the asymptotic direction for the solution $u(t)$.

In 2011, Zhang in \cite{16} proved the existence of the odd
symmetric parabolic or hyperbolic orbit for the restricted 3-body
problems with weak forces. He proved the following theorem.

\vspace{0.3cm}{\bf Theorem 1.3(See\cite{16})}\ \ {\em Suppose that
the potential $V(x)$ has the following form
\begin{eqnarray*}
V(x)=-\frac{1}{(|x|^{2}+r^{2})^{\kappa/2}},
\end{eqnarray*}
where $\kappa\in(0,2)$ and $r>0$. Then there exists one odd
parabolic or hyperbolic orbit for system $(1)-(2)$, which
minimizes the corresponding variational functional.}

In Theorem 1.3, the potential $V(x)$ has no singularity. When $V(x)$
is of class $C^{2}$, by taking the limit for a sequence of break
orbits, Serra \cite{11} obtained the existence of at least one
homoclinic orbit $u(t)$ at infinity, which means that
\begin{eqnarray}
\left\{
\begin{array}{ll}
\lim_{t\rightarrow\pm\infty}|u(t)|=+\infty\\
\lim_{t\rightarrow\pm\infty}\dot{u}(t)=0.
\end{array}
\right.\label{25}
\end{eqnarray}
Serra proved the following theorem.

\vspace{0.3cm}{\bf Theorem 1.4(See\cite{11})}\ \ {\em Suppose $V\in C^{2}(R^{N},R^{1})$ and satisfies\\

$(D_{1})$  $V(x)<0$ for all $x\in R^{N}$,

$(D_{2})$ there exist $R_{0}>0$ and $p>2$ such that

\begin{eqnarray*}
    V(x)=-\frac{1}{|x|^{p}}+W(x)\ \ \ \mbox{for all}\ \ \ |x|\geq R_{0},\ \ \mbox{with}\ W \mbox{ satisfying}
\end{eqnarray*}

$(D_{3})$ $\lim_{|x|\rightarrow +\infty} W(x)|x|^{p}=0$

$(D_{4})$ $(x,\nabla W(x))>0$ for all  $|x|\geq R_{0}$.\\

Then there exists at least one solution satisfying $(4)$
 for systems $(1)-(2)$ with $H=0$.}\\

Motivated by above papers, we study systems $(1)-(2)$, under some
weaker assumptions about the singularity for the potential, we
obtain the hyperbolic orbits
 with $H>0$. Precisely, we prove the
following theorem.

\vspace{0.3cm}{\bf Theorem 1.5}\ \ {\em Suppose that $V\in
C^{1}(R^{N}\setminus\{0\},R^{1})$ satisfies

\vspace{0.3cm}{\bf $(V_{1})$}\ $V(-x)=V(x)$, $\forall
x\in R^{N}\setminus\{0\}$,

\vspace{0.3cm}{\bf $(V_{2})$}\ $V(x)<0$, $\forall x\in
R^{N}\setminus\{0\}$,

\vspace{0.3cm}{\bf $(V_{3})$}\ $2 V(x)+(x,\nabla V(x))\rightarrow0\
\ \ \mbox{as}\ \ \ |x|\rightarrow+\infty$,

\vspace{0.3cm}{\bf $(V_{4})$}\ $2 V(x)+(x,\nabla
V(x))\rightarrow+\infty\ \ \ \mbox{as}\ \ \ |x|\rightarrow0$,

\vspace{0.3cm}{\bf $(V_{5})$}\ $-V(x)\rightarrow+\infty$ as
$|x|\rightarrow0$,

\vspace{0.3cm}{\bf $(V_{6})$}\ $V(x)\rightarrow0$ as
$|x|\rightarrow+\infty$.
\\

Then for any $H>0$, there is at least one hyperbolic orbit for
systems  $(1)-(2)$.}

Notice that the potential $V$ has no strong control at infinity even
though it satisfies $(V_{3})$ and $(V_{6})$. Under some additional
conditions, we can get the asymptotic direction of the solution at
infinity. We have the following theorem.

\vspace{0.3cm}{\bf Theorem 1.6}\ \ {\em Suppose that $V\in
C^{1}(R^{N}\setminus\{0\},R^{1})$ satisfies $(V_{1})$$-(V_{6})$ and
the following conditions

\vspace{0.3cm}{\bf $(V_{7})$}\ there exist constants $\beta>1$,
$M_{0}>0$ and $r_{0}\geq1$  such that
\begin{eqnarray*}
|x|^{\beta+1}|\nabla V(x)|\leq M_{0},\ \ \ |x|^{\beta+1}|V(x)|\leq
M_{0}\ \ \ \mbox{for all}\ \ \ |x|\geq r_{0}.
\end{eqnarray*}

Then for any $H>0$, there is at least one hyperbolic orbit for
systems  $(1)-(2)$ which has the given asymptotic direction at
infinity.}

\vspace{0.3cm}{\bf Remark 1}\ \ Notice that if
$V(x)=-\displaystyle\frac{1}{|x|^{\alpha}}$ $(\alpha>2)$, then $V$
satisfies $(V_{1})$$-$$(V_{6})$. It is easy to check that our
hypotheses are weaker than Theorem 1.1 and Theorem 1.2. Moreover,
there are functions which satisfy our hypotheses of Theorem 1.5 and
Theorem 1.6 but not $(A_{3})$. For example, let
\begin{eqnarray*}
V(x)=\left\{
\begin{array}{ll}
-\displaystyle\frac{\mbox{ln}(|x|^2+|x|^{-2})}{|x|^2}&\mbox{for $|x|\leq1,$}\\
J(x)&\mbox{for $1\leq|x|\leq2,$}\\
-\displaystyle\frac{1}{|x|^{3}} &\mbox{for $|x|\geq2,$}
\end{array}
\right.
\end{eqnarray*}
where $J(x)\in C^{1}(R^{N},R^{1})$ such that $V(x)\in
C^{1}(R^{N}\setminus\{0\},R^{1})$. The important difference between
our theorems and Theorem 1.3, Theorem 1.4 is that we have
singularities in Theorems 1.5 and 1.6.

\section{Variational Settings}

\ \ \ \ \ \ For any given unite vector(direction) $e\in S^{N-1}$, we
set
\begin{eqnarray*}
&&H^{1}=W^{1,2}(R^{1}/Z,R^{N}),\nonumber\\
&&E_{R}=\{q\in H^1|\ q(t+1/2)=-q(t), q(0)=q(1)=Re\},\nonumber\\
&&\Lambda_{R}=\{q\in E_{R}|\ q(t)\neq0,\forall t\in [0,1]\}.
\end{eqnarray*}

Here we just use $R$ to denote the Euclidean length of $q(0)$ and
$q(1)$. For any $q\in H^{1}$, we know that the following norms are
equivalent to each other
\begin{eqnarray}
&&\|q\|_{H^1}=\left(\int^{1}_{0}|\dot{q}(t)|^2dt\right)^{1/2}+\left|\int^{1}_{0}q(t)dt\right|\nonumber\\
&&\|q\|_{H^1}=\left(\int^{1}_{0}|\dot{q}(t)|^2dt\right)^{1/2}+\left(\int^{1}_{0}|q(t)|^2dt\right)^{1/2}\nonumber\\
&&\|q\|_{H^1}=\left(\int^{1}_{0}|\dot{q}(t)|^2dt\right)^{1/2}+|q(0)|.\nonumber
\end{eqnarray}
If $q\in \Lambda_{R}$, we have $\displaystyle\int^{1}_{0}q(t)dt=0$,
then by Poincar$\acute{\mbox{e}}$-Wirtinger's inequality, we obtain
that above norms are equivalent to
\begin{eqnarray*}
\|q\|_{H^1}=\left(\int^{1}_{0}|\dot{q}(t)|^2dt\right)^{1/2}.
\end{eqnarray*}

Let $L^{\infty}([0,1],R^{N})$ be a space of measurable functions from $[0,1]$ into $R^{N}$
and essentially bounded under the following norm
\begin{eqnarray*}
\ \|q\|_{L^{\infty}([0,1],R^{N})}:=esssup\{|q(t)|:t\in[0,1]\}.
\end{eqnarray*}
Moreover, let $f:\ \Lambda_{R}\rightarrow R^{1}$ be the functional
defined by
\begin{eqnarray*}
f(q)&=&\frac{1}{2}\int^{1}_{0}|\dot{q}(t)|^2dt\int^{1}_{0}(H-V(q(t)))dt\\
&=& \mbox{}\frac{1}{2}\|q\|^{2}\int^{1}_{0}(H-V(q(t)))dt
\end{eqnarray*}
Then one can easily check that $f\in C^{1}(\Lambda_{R},R^{1})$ and
\begin{eqnarray*}
\ \langle
f'(q),q(t)\rangle&=&\|q\|^2\int^{1}_{0}\left(H-V(q(t))-\frac{1}{2}(
\nabla V(q(t)),q(t))\right)dt.
\end{eqnarray*}

Our way to get the hyperbolic orbit is by approaching it with a
sequence of periodic solutions. Firstly, we prove the existence of
the approximate solutions, then we study the limit procedure.

\section{Existence of Periodic Solutions}

\ \ \ \ \ \ The approximate solutions are obtained by the
variational minimization methods. We need the following lemma which
is proved by A. Ambrosetti and V. Coti. Zelati in \cite{1}.

\vspace{0.3cm}{\bf Lemma 3.1(See\cite{1})}\ {\em Let
$f(q)=\displaystyle\frac{1}{2}\int^{1}_{0}|\dot{q}(t)|^2dt\int^{1}_{0}(H-V(q(t)))dt$
and $\tilde{q}\in H^{1}$ be such that $f^{'}(\tilde{q})=0$,
$f(\tilde{q})>0$. Set
\begin{eqnarray*}
T^{2}=\frac{\displaystyle\frac{1}{2}\displaystyle\int^{1}_{0}|\dot{\tilde{q}}(t)|^{2}dt}{\displaystyle\int^{1}_{0}(H-V(\tilde{q}(t))dt}.
\end{eqnarray*}
Then $\tilde{u}(t)=\tilde{q}(t/T)$ is a non-constant $T$-periodic
solution for (\ref{1}) and (\ref{2}).}

\vspace{0.3cm}{\bf Lemma 3.2(Palais\cite{13})}\ {\em Let $\sigma$ be
an orthogonal representation of a finite or compact group $G$ in the
real Hilbert space $H$ such that for any $\sigma\in G$,
\begin{eqnarray*}
f(\sigma\cdot x)=f(x),
\end{eqnarray*}
where $f\in C^{1}(H,R^{1})$. Let $S=\{x\in H|\sigma x=x,
\forall\sigma\in G\}$, then the critical point of $f$ in $S$ is also
a critical point of $f$ in $H$.}

\vspace{0.3cm}{\bf Lemma 3.3(Translation Property\cite{8})}\ {\em
Suppose that, in domain $D\subset R^{N}$, we have a solution
$\phi(t)$ for the following differential equation
\begin{eqnarray*}
x^{(n)}+F(x^{(n-1)},\cdots,x)=0,
\end{eqnarray*}
where $x^{(k)}=d^{k}x/dt^{k}$, $k=0,1,\cdots,n$, $x^{(0)}=x$. Then
$\phi(t-t_{0})$ with $t_{0}$ being a constant is also a solution.}

In the following, we introduce $Gordon's\ Strong\ Force$ condition.

\vspace{0.3cm}{\bf Lemma 3.4(Gordon\cite{19})}\ \ {\em $V$ is said
to satisfy $Gordon's\ Strong\ Force$ condition, if there exists a
neighborhood $\mathcal{N}$ of 0 and a function $U\in
C^{1}(R^{N}\setminus\{0\},R^{1})$ such that}

$(\mbox{i})$\ $\lim_{x\to 0}U(x)=-\infty$;

$(\mbox{ii})$\ $-V(x)\geq |U'(x)|^{2}$ for every $x\in \mathcal{N}$\
$\setminus\{0\}$.
\\

If $V$ satisfies $Gordon's\ Strong\ Force$ condition, then
\begin{eqnarray*}
\int^{1}_{0}V(x_{j})dt\rightarrow-\infty,\ \ \ \forall\
x_{j}\rightharpoonup x\in \partial \Lambda_{R}.
\end{eqnarray*}

\vspace{0.3cm}{\bf Lemma 3.5} {\em Suppose $(V_2)$ and $(V_4)$ hold,
then $V$ satisfies $Gordon's\ Strong\ Force$ condition.}

 \vspace{0.3cm}{\bf Proof.}\ Let $\phi(r)=-V(r\widetilde{e})r^{2}$, where $r=|x|$, $\widetilde{e}=x/|x|$, then we have
\begin{eqnarray*}
\phi'(r)=-r(2V(r\widetilde{e})+(\nabla
V(r\widetilde{e}),r\widetilde{e})).
\end{eqnarray*}
It follows from $(V_{4})$ that, there exists a constant $\delta>0$
such that
\begin{eqnarray*}
\phi'(r)\leq0\ \ \ \mbox{for all}\ \ \ 0<|x|\leq\delta.
\end{eqnarray*}
Since $V\in C^{1}(R^{N}\setminus\{0\},R^{1})$, we get
\begin{eqnarray*}
\phi(r)&\geq&\phi(\delta)=-V(\delta e)\delta^{2}\geq
\delta^{2}\min_{|x|=\delta}(-V(x)).
\end{eqnarray*}
It follows from the definition of $\phi$ and $(V_2)$ that there
exists a constant $C>0$ such that
\begin{eqnarray*}
-V(x)\geq\frac{C}{|x|^{2}}\ \ \ \mbox{for all}\ \ \ 0<r\leq\delta.
\end{eqnarray*}
We set $U(x)=\sqrt{C}\ln|x|$, then by some calculation, we obtain
\begin{eqnarray*}
\lim_{x\to 0}U(x)=-\infty\ \ \ \mbox{and}\ \ \ -V(x)\geq|U'(x)|^{2}\
\ \ \mbox{for all}\ \ \ 0<r\leq\delta,
\end{eqnarray*}
which proves this lemma.

\vspace{0.3cm}{\bf Lemma 3.6}\ {\em Suppose the conditions of
Theorem 1.5 hold, then for any  $R>0$, there exists at least one
periodic solution on $\Lambda_{R}$ for the following systems
\begin{equation}
   \ddot{u}(t)+\nabla V(u(t))=0,\ \ \ \
   \forall\ t\in\left(-\frac{T_{R}}{2},\frac{T_{R}}{2}\right)\label{16}
\end{equation}
with
\begin{equation}
   \frac{1}{2}|\dot{u}(t)|^{2}+ V(u(t))=H,\ \ \ \ \ \
   \forall\ t\in\left(-\frac{T_{R}}{2},\frac{T_{R}}{2}\right).\label{17}
\end{equation}}

\vspace{0.3cm}{\bf Proof.}\ We notice that $H^{1}$ is a reflexive
Banach space and $E_{R}$ is a weakly closed subset of $H^{1}$. Since
$H>0$, we obtain that
\begin{eqnarray}
f(q)&=&\frac{1}{2}\|q\|^{2}\int^{1}_{0}(H-V(q(t))dt\geq\frac{H}{2}\|q\|^{2},\label{9}
\end{eqnarray}
which implies that $f$ is a functional bounded from below,
furthermore, it is easy to check that $f$ is weakly lower
semi-continuous and
\begin{eqnarray}
f(q)\rightarrow +\infty\ \ \ \mbox{as}\ \ \|q\|\rightarrow
+\infty.
\end{eqnarray}

Then, we conclude that for every $R>0$ there exists a minimizer
$q_{R}\in E_{R}$ such that
\begin{eqnarray}
f'(q_R)=0,\ \ \ \ f(q_R)=\inf_{q\in E_{R}} f(q)>0.\label{18}
\end{eqnarray}
Furthermore, we need to prove that $q_{R}\in\Lambda_{R}$ which
means $q_{R}$ has no collision for any $R>0$. Suppose that
$\{q_{j}\}_{j\in N}$ is the minimizing sequence, then if $q_{R}$
has collision, which means $q_{R}\in
\partial\Lambda_{R}=\{q_{R}\in E_{R}|\ \exists\ t'\in [0,1]\ st.\
q_{R}(t')=0\}$, we can prove that
\begin{eqnarray}
f(q_{j})\rightarrow +\infty,\ \ \ \mbox{as}\ \
j\rightarrow+\infty.\label{10}
\end{eqnarray}
To prove this fact, there are two cases needed to be discussed.

{\bf Case 1.}\ If $q_R=$constant, it follows from $q_{R}\in
\partial\Lambda_{R}$ that $q_{R}\equiv0$, which is a contradiction,
since $|q_{R}(0)|=|q_{R}(1)|=R$.

{\bf Case 2.}\ If $q_{R}\neq$constant, we have
$\|q_{R}\|^{2}=\displaystyle\int^{1}_{0}|\dot{q}_{R}(t)|^{2}dt>0$,
otherwise by $q_{R}(t+1/2)=-q_{R}(t)$, we can deduce
$q_{R}\equiv0$ which is a contradiction. Then by the
weakly-lower-semi-continuity of norm, we have
\begin{eqnarray*}
\liminf_{j\rightarrow\infty}\|q_{j}\|\geq \|q_{R}\|>0.
\end{eqnarray*}
Then by Lemma 3.4, (\ref{10}) holds.

Moreover, let $Q(t)=R(\xi\cos(2\pi t)+\eta\sin(2\pi t))\in
\Lambda_R$, where $\xi, \eta \in R^{N}\setminus\{0\}$,
$|\xi|=|\eta|=1$, $(\xi,\eta)=0$,  which implies that $|Q(t)|=R$,
$\|Q\|=2\pi R$, hence
\begin{eqnarray*}
f(Q)&=&2\pi^{2}R^{2}\left(H-\int^{1}_{0}V(Q(t))dt\right).
\end{eqnarray*}
Since $V\in C^{1}(R^{N}\setminus\{0\},R^{1})$, then there exists a
constant $M_{1,R}> 0$ such that $|V(Q(t))|\leq M_{1,R}$. We obtain
that
\begin{eqnarray}
f(q_R)\leq f(Q)\leq M_{2,R}\label{22}
\end{eqnarray}
for some $M_{2,R}>0$, but (\ref{22}) contradicts with (\ref{10}) for
any fixed $R>0$. Then we can see that $q_{R}\in \Lambda_{R}$ has no
collision.

Let
\begin{eqnarray*}
T_{R}^{2}=\frac{\displaystyle\frac{1}{2}\displaystyle\int^{1}_{0}|\dot{q}_{R}(t)|^{2}dt}{\displaystyle\int^{1}_{0}(H-V(q_{R}(t))dt}.
\end{eqnarray*}

Then by Lemma 3.1$-$ Lemma 3.3, we obtain that
$u_{R}(t)=q_{R}(\frac{t+\frac{T_{R}}{2}}{T_{R}}):\left(-\frac{T_{R}}{2},\frac{T_{R}}{2}\right)\rightarrow
\Lambda_{R}$ is a $T_{R}$-periodic solution for systems (\ref{16})
and (\ref{17}). The lemma is proved.

\section{Blowing-up Argument}

\ \ \ \ \ \ Subsequently, we need to show that $u_{R}(t)$ can not
diverge to infinity uniformly as $R\rightarrow+\infty$. Moreover, we
prove the following lemma.

\vspace{0.3cm}{\bf Lemma 4.1}\ {\em Suppose that
$u_{R}(t):\left(-\frac{T_{R}}{2},\frac{T_{R}}{2}\right)\rightarrow
\Lambda_{R}$ is the solution obtained in Lemma 3.6, then
$\min_{t\in\left[-\frac{T_{R}}{2},\frac{T_{R}}{2}\right]}|u_{R}(t)|$
is bounded from above. More precisely, there is a constant $M>0$
independent of $R$ such that
\begin{eqnarray*}
\min_{t\in\left[-\frac{T_{R}}{2},\frac{T_{R}}{2}\right]}|u_{R}(t)|\leq
M\ \ \ \mbox{for all}\ \ \ R>0.
\end{eqnarray*} }

\vspace{0.3cm}{\bf Proof.}\ Since $q_{R}\in \Lambda_{R}$, it is easy
to see that $u_{R}(t)=q_{R}(\frac{t+\frac{T_{R}}{2}}{T_{R}})$
satisfies $u_{R}(-\frac{T_{R}}{2})=u_{R}(\frac{T_{R}}{2})$ and
$\dot{u}_{R}(-\frac{T_{R}}{2})=\dot{u}_{R}(\frac{T_{R}}{2})$, then
we have that
\begin{eqnarray*}
&&\left(u_{R}\left(\frac{T_{R}}{2}\right),\dot{u}_{R}\left(\frac{T_{R}}{2}\right)\right)
-\left(u_{R}\left(-\frac{T_{R}}{2}\right),\dot{u}_{R}\left(-\frac{T_{R}}{2}\right)\right)\nonumber\\
&=& \mbox {}
\int^{\frac{T_{R}}{2}}_{-\frac{T_{R}}{2}}\frac{d}{dt}(u_{R}(t),\dot{u}_{R}(t))dt\nonumber\\
&=& \mbox {} \int^{\frac{T_{R}}{2}}_{-\frac{T_{R}}{2}}
(|\dot{u}_{R}(t)|^{2}+(u_{R}(t),\ddot{u}_{R}(t)))dt\nonumber\\
&=& \mbox {}
\int^{\frac{T_{R}}{2}}_{-\frac{T_{R}}{2}}2(H-V(u_{R}(t)))-(\nabla
V(u_{R}(t)),u_{R}(t))dt.
\end{eqnarray*}
Then we obtain that
\begin{eqnarray*}
\int^{\frac{T_{R}}{2}}_{-\frac{T_{R}}{2}}2H-(2V(u_{R}(t))+(\nabla
V(u_{R}(t)),u_{R}(t)))dt=0.
\end{eqnarray*}
There are two cases needed to be discussed.

{\bf Case 1.}\ $2H-(2V(u_{R}(t))+(\nabla
V(u_{R}(t)),u_{R}(t)))\equiv0$, which implies that
\begin{eqnarray*}
2H&=&2V(u_{R}(t))+(\nabla V(u_{R}(t)),u_{R}(t)),\ \ \ \mbox{a.e.} \
t\in\left[-\frac{T_{R}}{2},\frac{T_{R}}{2}\right].
\end{eqnarray*}
Hypothesis $(V_{3})$ implies that there exists a constant
$M_{1}>0$ independent of $R$ such that
\begin{eqnarray*}
\min_{t\in\left[-\frac{T_{R}}{2},\frac{T_{R}}{2}\right]}|u_{R}(t)|\leq
M_{1}.
\end{eqnarray*}

{\bf Case 2.}\ $2(H-V(u_{R}(t)))-(\nabla V(u_{R}(t)),u_{R}(t))$
changes sign in $\left[-\frac{T_{R}}{2},\frac{T_{R}}{2}\right]$.
Then there exists
$t_{0}\in\left[-\frac{T_{R}}{2},\frac{T_{R}}{2}\right]$ such that
\begin{eqnarray*}
2H-(2V(u_{R}(t_{0}))+(\nabla V(u_{R}(t_{0})),u_{R}(t_{0})))<0,
\end{eqnarray*}
which implies that
\begin{eqnarray*}
2H&<&2V(u_{R}(t_{0}))+(\nabla V(u_{R}(t_{0})),u_{R}(t_{0})).
\end{eqnarray*}
It follows from hypothesis $(V_{3})$ that there exists a constant
$M_{2}>0$ independent of $R$ such that
\begin{eqnarray*}
\min_{t\in\left[-\frac{T_{R}}{2},\frac{T_{R}}{2}\right]}|u_{R}(t)|\leq
M_{2}.
\end{eqnarray*}
Then the proof is completed.

\section{Proof of Theorem 1.5}

\ \ \ \ \ \ The ideas for the following proofs in this section
mostly comes from Lemma 2.1 and Lemma 4.1 in \cite{3}, we write out
them for completeness.

\vspace{0.3cm}{\bf Lemma 5.1}\ {\em Suppose that $u_{R}(t)$ is the
solution for $(\ref{16})-(\ref{17})$ obtained in Lemma 3.6. Then
there exists a constant $m>0$ independent of $R$ such that}
\begin{eqnarray*}
\min_{t\in\left[-\frac{T_{R}}{2},\frac{T_{R}}{2}\right]}|u_{R}(t)|\geq
m.
\end{eqnarray*}

\vspace{0.3cm}{\bf Proof.}\ Since $u_{R}(t)$ is a solution for
system $(\ref{16})-(\ref{17})$, then we can deduce that
\begin{eqnarray*}
\frac{d^{2}}{dt^{2}}\frac{1}{2}|u_{R}(t)|^{2}&=&\frac{d}{dt}(u_{R}(t),\dot{u}_{R}(t))\nonumber\\
&=& \mbox{} |\dot{u}_{R}(t)|^{2}+(u_{R}(t),\ddot{u}_{R}(t))\nonumber\\
&=& \mbox{} 2(H-V(u_{R}(t)))-(\nabla V(u_{R}(t)),u_{R}(t))\ \ \ \
t\in\left(-\frac{T_{R}}{2},\frac{T_{R}}{2}\right).
\end{eqnarray*}
Since
$\left|u_{R}\left(-\frac{T_{R}}{2}\right)\right|=\left|u_{R}\left(\frac{T_{R}}{2}\right)\right|=R
$, then using hypothesis $(V_{4})$, we can find $m\in (0,1)$
independent of $R$ such that, for any
$t\in\{t\in\left[-\frac{T_{R}}{2},\frac{T_{R}}{2}\right]
|\max_{t\in\left[-\frac{T_{R}}{2},\frac{T_{R}}{2}\right]}|u_{R}(t)|\leq
m\}$,
\begin{eqnarray*}
2H-(2V(u_{R}(t))+(\nabla V(u_{R}(t)),u_{R}(t)))<0,
\end{eqnarray*}
which implies that $|u_{R}(t)|$ is concave when $|u_{R}(t)|\leq m$
and $|u_{R}(t)|$ cannot take a local minimum such that
$\max_{t\in\left[-\frac{T_{R}}{2},\frac{T_{R}}{2}\right]}|u_{R}(t)|\leq
m$, which implies that
\begin{eqnarray*}
 |u_{R}(t)|\geq m\ \ \ \mbox{for}\ \ \mbox{all}\ \ \
 t\in\left[-\frac{T_{R}}{2},\frac{T_{R}}{2}\right].
\end{eqnarray*}
If not, we can assume that there exists a
$\overline{t}\in\left[-\frac{T_{R}}{2},\frac{T_{R}}{2}\right]$ such
that $|u_{R}(\overline{t})|<m$, then we can easily check that
$|u_{R}(t)|$ takes a local minimum at some $\tilde{t}$ with
$|u_{R}(\tilde{t})|<m$, which is a contradiction. Then we obtain the
conclusion.

\vspace{0.3cm}{\bf Lemma 5.2}\ {\em Suppose that $R>M$, where $M$ is
defined in Lemma 4.1 and $u_{R}(t)$ is the solution for
$(\ref{16})-(\ref{17})$ obtained in Lemma 3.6. Set
\begin{eqnarray*}
t_{+}=\sup\left\{t\in\left[-\frac{T_{R}}{2},\frac{T_{R}}{2}\right]|\left|u_{R}(t)\right|\leq
L\right\}
\end{eqnarray*}
and
\begin{eqnarray*}t_{-}=\inf\left\{t\in\left[-\frac{T_{R}}{2},\frac{T_{R}}{2}\right]|\left|u_{R}(t)\right|\leq
L\right\}
\end{eqnarray*}
where $L$ is a constant independent of $R$ such that $M<L<R$. Then
we have that}
\begin{eqnarray*}
\frac{T_{R}}{2}-t_{+}\rightarrow+\infty,\ \ \
t_{-}+\frac{T_{R}}{2}\rightarrow+\infty\ \ \ \mbox{as}\ \
R\rightarrow+\infty.
\end{eqnarray*}

\vspace{0.3cm}{\bf Proof.}\ By the definition of $u_{R}(t)$ we have
that
\begin{eqnarray*}
\left|u_{R}\left(-\frac{T_{R}}{2}\right)\right|=\left|u_{R}\left(\frac{T_{R}}{2}\right)\right|=R.
\end{eqnarray*}
Then, by $(V_{2})$ and the definitions of $t_{+}$ and $t_{-}$, we
have
\begin{eqnarray}
\int^{\frac{T_{R}}{2}}_{t_{+}}\sqrt{H-V(u_{R}(t))}|\dot{u}_{R}(t)|dt&\geq&\sqrt{H}\int^{\frac{T_{R}}{2}}_{t_{+}}|\dot{u}_{R}(t)|dt\nonumber\\
&\geq& \mbox{}
\sqrt{H}\left|\int^{\frac{T_{R}}{2}}_{t_{+}}\dot{u}_{R}(t)dt\right|\geq\sqrt{H}(R-L)\label{3}
\end{eqnarray}
and
\begin{eqnarray}
\int^{t_{-}}_{-\frac{T_{R}}{2}}\sqrt{H-V(u_{R}(t))}|\dot{u}_{R}(t)|dt&\geq&\sqrt{H}\int^{t_{-}}_{-\frac{T_{R}}{2}}|\dot{u}_{R}(t)|dt
\nonumber\\
&\geq& \mbox{}
\sqrt{H}\left|\int^{t_{-}}_{-\frac{T_{R}}{2}}\dot{u}_{R}(t)dt\right|\geq\sqrt{H}(R-L)\label{4}.
\end{eqnarray}
Since $V\in C^{1}(R^{N}\setminus\{0\},R^{1})$, it follows from Lemma
5.1 and $(V_{6})$, that there exists a constant $M_{3}>0$
independent of $R$ such that
\begin{eqnarray*}
|V(u_{R}(t))|\leq M_{3}\ \ \ \mbox{for}\ \ \mbox{all}\ \
t\in\left[-\frac{T_{R}}{2},\frac{T_{R}}{2}\right],
\end{eqnarray*}
which implies that
\begin{eqnarray*}
\int^{\frac{T_{R}}{2}}_{t_{+}}\sqrt{H-V(u_{R}(t))}|\dot{u}_{R}(t)|dt=\sqrt{2}\int^{\frac{T_{R}}{2}}_{t_{+}}(H-V(u_{R}(t)))dt\leq
\sqrt{2}(H+M_{3})\left(\frac{T_{R}}{2}-t_{+}\right).
\end{eqnarray*}
Combining (\ref{3}) with the  above estimates, we obtain that
\begin{eqnarray*}
\sqrt{H}(R-L)\leq\sqrt{2}(H+M_{3})\left(\frac{T_{R}}{2}-t_{+}\right).
\end{eqnarray*}
Then we have
\begin{eqnarray*}
\frac{T_{R}}{2}-t_{+}\rightarrow+\infty,\ \ \ \mbox{as}\ \
R\rightarrow+\infty.
\end{eqnarray*}
The limit for $t_{-}+\frac{T_{R}}{2}$ can be obtained in the similar
way. The proof is completed.

Subsequently, we set that
\begin{eqnarray}
t^{*}=\inf\left\{t\in\left[-\frac{T_{R}}{2},\frac{T_{R}}{2}\right]||u_{R}(t)|=M\right\}\label{31}
\end{eqnarray}
and
\begin{eqnarray*}
u_{R}^{*}(t)=u_{R}(t-t^{*}).
\end{eqnarray*}
Since $L>M$, we can deduce that $t_{+}\geq t^{*}\geq t_{-}$, which
implies that
\begin{eqnarray*}
-\frac{T_{R}}{2}+t^{*}\rightarrow-\infty,\ \ \mbox{}\ \
\frac{T_{R}}{2}+t^{*}\rightarrow+\infty\ \ \mbox{as}\ \
R\rightarrow\infty.
\end{eqnarray*}

Then it follows from $(\ref{17})$ that
\begin{eqnarray*}
   \frac{1}{2}|\dot{u}^{*}_{R}(t)|^{2}+ V(u^{*}_{R}(t))=H,\ \ \
   \forall\
   t\in\left(-\frac{T_{R}}{2}+t^{*},\frac{T_{R}}{2}+t^{*}\right),
\end{eqnarray*}
which implies that
\begin{eqnarray*}
|\dot{u}^{*}_{R}(t)|^{2}=2(H-V(u^{*}_{R}(t))),\ \ \
  \forall\ t\in\left(-\frac{T_{R}}{2}+t^{*},\frac{T_{R}}{2}+t^{*}\right).
\end{eqnarray*}
By Lemma 5.1,  $(V_{6})$ and $V\in
C^{1}(R^{N}\setminus\{0\},R^{1})$, we can deduce that there exists
a constant $M_{4}>0$ independent of $R$ such that
\begin{eqnarray*}
 |V(u^{*}_{R}(t))|\leq M_{4}\ \ \ \mbox{for all}\ \
 t\in\left(-\frac{T_{R}}{2}+t^{*},\frac{T_{R}}{2}+t^{*}\right).
\end{eqnarray*}
Then there is a constant $M_{5}>0$ independent of $R$ such that
\begin{eqnarray*}
 |\dot{u}^{*}_{R}(t)|\leq M_{5}\ \ \ \mbox{for all}\ \
 t\in\left(-\frac{T_{R}}{2}+t^{*},\frac{T_{R}}{2}+t^{*}\right).
\end{eqnarray*}
which implies that
\begin{eqnarray}
\ |u_{R}^{*}(t_{1})-u_{R}^{*}(t_{2})| \leq
\left|\int^{t_{1}}_{t_{2}}\dot{u}_{R}^{*}(s)ds\right| \leq \mbox{}
\int^{t_{1}}_{t_{2}}\left|\dot{u}_{R}^{*}(s)\right|ds \leq
M_{5}|t_{1}-t_{2}|\label{27}
\end{eqnarray}
for each $R>0$ and $t_{1}, t_{2} \in
\left(-\frac{T_{R}}{2}+t^{*},\frac{T_{R}}{2}+t^{*}\right)$, which
shows $\{u_{R}^{*}\}$ is equicontinuous. Then there is a subsequence
$\{u_{R}^{*}\}_{R>0}$ converging to $u_{\infty}$ in
$C_{loc}(R^{1},R^{N})$. Then there exists a function $u_{\infty}(t)$
such that
\begin{eqnarray*}
&&(\mbox{i})\ u_{R}^{*}(t)\rightarrow u_{\infty}(t)\ \mbox{in}\ C_{loc}(R^{1},R^{N})\nonumber\\
&&(\mbox{ii})|u_{\infty}(t)|\rightarrow +\infty\ \mbox{as}\
|t|\rightarrow+\infty
\end{eqnarray*}
and $u_{\infty}(t)$ satisfies systems $(1)-(2)$. Then we finish the
proof of Theorem 1.5.

\section{Proof of Theorem 1.6}
\ \ \ \ \ \ By the conditions of Theorem 1.6, the existence of
hyperbolic solutions for systems $(1)-(2)$ can be obtained with a
similar proof of Theorem 1.5. Subsequently, we give the proof of the
asymptotic direction of hyperbolic solutions at infinity. The proof
is similar to Felmer and Tanaka's in \cite{3}.

\vspace{0.3cm}{\bf Lemma 6.1}\ {\em Suppose that $u_{R}(t)$ is the
solution for $(\ref{16})-(\ref{17})$ obtained in Lemma 3.6. Then
there exists a constant $M_{6}>0$ independent of $R>1$ such that}
\begin{eqnarray*}
\int^{\frac{T_{R}}{2}}_{-\frac{T_{R}}{2}}\sqrt{H-V(u_{R}(t))}|\dot{u}_{R}(t)|dt\leq
\sqrt{2H}R+M_{6}.
\end{eqnarray*}

\vspace{0.3cm}{\bf Proof.}\ Firstly, we define the function $\xi(t)$
on $[1,+\infty)$ as a solution of
\begin{eqnarray*}
&& \dot{\xi}(t)=\sqrt{2(H-V(\xi(t)e))}\\
&& \xi(1)=1.
\end{eqnarray*}
And $\tau_{R}>1$ is a real number such that $\xi(\tau_{R})=R$. We
can define $\xi(t)$ in $(-\infty,0]$ and $\tau_{-R}$ in a similar
way. Then we can fix $\varphi(t)\in H^{1}([0,1],R^{N})$ such that
$\tilde{\gamma}_{R}(t)\in \Lambda_{R}$ where
\begin{eqnarray*}
\tilde{\gamma}_{R}(t)=\gamma_{R}(t(\tau_{R}-\tau_{-R})+\tau_{-R}),\
\ \ \mbox{and} \ \ \gamma_{R}(t)=\left\{
\begin{array}{ll}
\xi(t)e&\mbox{for $t\in[1,\tau_{R}]\bigcup[\tau_{-R},0]$}\\
\varphi(t)&\mbox{for $t\in[0,1]$}.
\end{array}
\right.
\end{eqnarray*}
Subsequently, we set $u_{r}(t)=\tilde{\gamma}_{R}(\frac{t+r}{2r})$.
And it is easy to see that $u_{r}(t)=\gamma_{R}(t)$ if $\pm
r=\tau_{\pm R}$. Similar to \cite{3}, we can deduce that for $r>0$
\begin{eqnarray}
(2f(\tilde{\gamma}_{R}))^{\frac{1}{2}}&=&\inf_{r>0}\frac{1}{\sqrt{2}}\int^{r}_{-r}\frac{1}{2}|\dot{u}_{r}(t)|^{2}+H-V(u_{r}(t))dt\nonumber\\
&\leq& \mbox{}
\frac{1}{\sqrt{2}}\int^{\tau_{R}}_{-\tau_{R}}\frac{1}{2}|\dot{\gamma}_{R}(t)|^{2}+H-V(\gamma_{R}(t))dt.\label{26}
\end{eqnarray}
Since
$[-\tau_{R},\tau_{R}]=[-\tau_{R},0]\bigcup[0,1]\bigcup[1,\tau_{R}]$,
then by $(V_{7})$, we can estimate (\ref{26}) by three integral.
Firstly, we estimate the integral on $[1,\tau_{R}]$, which is
\begin{eqnarray*}
I_{[1,\tau_{R}]}&=&\frac{1}{\sqrt{2}}\int^{\tau_{R}}_{1}\frac{1}{2}|\dot{\gamma}_{R}(t)|^{2}+H-V(\gamma_{R}(t))dt\nonumber\\
&=& \mbox{}\int^{\tau_{R}}_{1}\sqrt{H-V(\xi(t)e)}\dot{\xi}(t)dt\nonumber\\
&=& \mbox{} \int^{R}_{1} \sqrt{H-V(se)}ds\nonumber\\
&\leq& \mbox{} \int^{R}_{1} \sqrt{H}+\sqrt{-V(se)}ds\nonumber\\
&=& \mbox{}
\sqrt{H}(R-1)+\int^{r_{0}}_{1}\sqrt{-V(se)}ds+\int^{R}_{r_{0}}\sqrt{-V(se)}ds\nonumber\\
&\leq& \mbox{}
\sqrt{H}(R-1)+M_{4}(r_{0}-1)+\sqrt{M_{0}}\int^{R}_{r_{0}}s^{-\frac{\beta+1}{2}}ds\nonumber\\
&\leq& \mbox{}\sqrt{H}(R-1)+M_{4}(r_{0}-1)+\sqrt{M_{0}}\int^{+\infty}_{r_{0}}s^{-\frac{\beta+1}{2}}ds\nonumber\\
&\leq& \mbox{} \sqrt{H}R + M_{7}
\end{eqnarray*}
for some $M_{7}>0$ independent of $R\geq1$. Similarly, we can get
\begin{eqnarray*}
I_{[\tau_{R},0]}\leq \sqrt{H}R + M_{7}.
\end{eqnarray*}
Since $I_{[0,1]}$ is independent of $R$, we obtain that
\begin{eqnarray*}
\frac{1}{\sqrt{2}}\int^{\tau_{R}}_{-\tau_{R}}\frac{1}{2}|\dot{\gamma}_{R}(t)|^{2}+H-V(\gamma_{R}(t))dt\leq
2\sqrt{H}R + M_{6}
\end{eqnarray*}
for some $M_{6}>0$ independent of $R$. Then by (\ref{26}) and $q(t)$
is the minimizer of $f$ on $\Lambda_{R}$, we have
\begin{eqnarray}
\int^{\frac{T_{R}}{2}}_{-\frac{T_{R}}{2}}\sqrt{H-V(u_{R}(t))}|\dot{u}_{R}(t)|dt
&\leq&
\left(\int^{\frac{T_{R}}{2}}_{-\frac{T_{R}}{2}}H-V(u_{R}(t))dt\right)^{\frac{1}{2}}\left(\int^{\frac{T_{R}}{2}}_{-\frac{T_{R}}{2}}|\dot{u}_{R}(t)|^{2}dt\right)^{\frac{1}{2}}
\nonumber\\
&=&\mbox{}(2f(q))^{\frac{1}{2}}\nonumber\\
&\leq&\mbox{}(2f(\tilde{\gamma}_{R}))^{\frac{1}{2}}\nonumber\\
&\leq& \mbox{}
\frac{1}{\sqrt{2}}\int^{\tau_{R}}_{-\tau_{R}}\frac{1}{2}|\dot{\gamma}_{R}(t)|^{2}+H-V(\gamma_{R}(t))dt\nonumber\\
&\leq& \mbox{}2\sqrt{H}R + M_{6}.
\end{eqnarray}
Then we finish the proof of this lemma.

Similar to Felmer and Tanaka \cite{3}, we set
\begin{eqnarray*}
A(t)=\sqrt{|u_{R}(t)|^{2}|\dot{u}_{R}(t)|^{2}-(u_{R}(t),\dot{u}_{R}(t))^{2}}
\end{eqnarray*}
and
\begin{eqnarray*}
\omega(t)=\frac{A(t)}{|u_{R}(t)||\dot{u}_{R}(t)|}.
\end{eqnarray*}
Using the motion and energy equations, we have
\begin{eqnarray*}
|\dot{A}(t)|\leq |u_{R}(t)||\nabla V(u_{R}(t))|
\end{eqnarray*}
and
\begin{eqnarray*}
\frac{d\omega}{dt}=\frac{2}{|u_{R}(t)||\dot{u}_{R}(t)|}(-\omega\sqrt{1-\omega^{2}}sign(u_{R}(t),\dot{u}_{R}(t))(H-V(u_{R}(t)))+|u_{R}(t)||\nabla
V(u_{R}(t))|).
\end{eqnarray*}
The proof of the following lemma is the same as \cite{3}.

\vspace{0.3cm}{\bf Lemma 6.2(See\cite{3})}\ {\em Assume $u_{R}$ is a
solution for $(\ref{16})-(\ref{17})$ obtained in Lemma 3.6. For any
$\eta\in(0,1)$, there exists a $L_{\eta}\geq m$ such that if}
\begin{eqnarray}
|u_{R}(t_{0})|\geq L_{\eta},\ \ \
(u_{R}(t_{0}),\dot{u}_{R}(t_{0}))>0\ \ \ \mbox{and}\ \ \
\omega(t_{0})<\eta\label{28}
\end{eqnarray}
{\em for some $t_{0}\in (-\frac{T_{R}}{2},\frac{T_{R}}{2})$, then we
have for $t\in[t_{0},\frac{T_{R}}{2}]$}
\begin{eqnarray*}
&&(\mbox{i}).\ \omega(t)<\eta,\nonumber\\
&&(\mbox{ii}).\ \frac{d}{dt}|u_{R}(t)|\geq\sqrt{1-\eta^{2}}|\dot{u}_{R}(t)|,\nonumber\\
&&(\mbox{iii}).\ \frac{d}{dt}|u_{R}(t)|\geq\sqrt{2(1-\eta^{2})H},\nonumber\\
&&(\mbox{iv}).\
|u_{R}(t)|\geq|u_{R}(t_{0})|+\sqrt{2(1-\eta^{2})H}(t-t_{0}).
\end{eqnarray*}

\vspace{0.3cm}{\bf Lemma 6.3(See\cite{3})}\ {\em Let $u_{R}$ is a
solution for $(\ref{16})-(\ref{17})$ obtained in Lemma 3.6
satisfying (\ref{28}) and $|u_{R}(t)|\geq r_{0}$ with $t\geq t_{0}$
for certain $t_{0}\in (-\frac{T_{R}}{2},\frac{T_{R}}{2})$ with
$\eta\in(0,\frac{1}{2})$ and $L_{\eta}$ as in Lemma 6.2. Then for
$t\geq t_{0}$ we have}
\begin{eqnarray*}
\left|\frac{u_{R}(t)}{|u_{R}(t)|}-\frac{u_{R}(t_{0})}{|u_{R}(t_{0})|}\right|\leq
M_{8}\eta+\frac{M_{9}}{|u_{R}(t_{0})|^{\beta}},
\end{eqnarray*}
where $M_{8}$, $M_{9}>0$ are independent of $\eta$, $u_{R}(t)$ and
$t_{0}$.

\vspace{0.3cm}{\bf Proof.}\ By Lemma 5.1, (iii) of Lemma 6.2 and
$(V_{7})$, we can estimate $A(t)$ as following.
\begin{eqnarray}
A(t)&\leq& A(t_{0})+\int^{t}_{t_{0}}|u_{R}(s)||\nabla V(u_{R}(s))|ds\nonumber\\
&\leq& \mbox{}
A(t_{0})+\frac{M_{9}}{\sqrt{2(1-\eta^{2})H}}\int^{t}_{t_{0}}|u_{R}(s)||\nabla
V(u_{R}(s))|\frac{d}{ds}|u_{R}(s)|ds\nonumber\\
&\leq& \mbox{}
A(t_{0})+\frac{M_{9}}{\sqrt{2(1-\eta^{2})H}}\int^{|u_{R}(t)|}_{|u_{R}(t_{0})|}\varphi\left|\nabla
V\left(\varphi \frac{u_{R}(s)}{|u_{R}(s)|}\right)\right|d\varphi\nonumber\\
&\leq& \mbox{}
A(t_{0})+\frac{M_{9}M_{0}}{\sqrt{2(1-\eta^{2})H}}\int^{|u_{R}(t)|}_{|u_{R}(t_{0})|}\frac{1}{\varphi^{\beta}}d\varphi\nonumber\\
&\leq& \mbox{}
A(t_{0})+\frac{M_{9}M_{0}}{\sqrt{2(1-\eta^{2})H}(\beta-1)}\frac{1}{|u_{R}(t_{0})|^{\beta-1}}\nonumber\\
&\leq& \mbox{}A(t_{0})+
\frac{M_{10}}{|u_{R}(t_{0})|^{\beta-1}}\label{29}
\end{eqnarray}
for some $M_{10}>0$ independent of $R$. Since we have
\begin{eqnarray}
\left|\frac{d}{dt}\frac{u_{R}(t)}{|u_{R}(t)|}\right|=\frac{A(t)}{|u_{R}(t)|^{2}},\label{30}
\end{eqnarray}
then it follows from (iii) of Lemma 6.2, (\ref{29}) and (\ref{30})
that
\begin{eqnarray*}
\left|\frac{u_{R}(t)}{|u_{R}(t)|}-\frac{u_{R}(t_{0})}{|u_{R}(t_{0})|}\right|&\leq&\int^{t}_{t_{0}}\frac{A(s)}{|u_{R}(s)|^{2}}ds\nonumber\\
&\leq& \mbox{}\left(A(t_{0})+
\frac{M_{10}}{|u_{R}(t_{0})|^{\beta-1}}\right)\int^{t}_{t_{0}}\frac{1}{|u_{R}(s)|^{2}}ds\nonumber\\
&\leq& \mbox{}\left(A(t_{0})+
\frac{M_{10}}{|u_{R}(t_{0})|^{\beta-1}}\right)\frac{1}{\sqrt{2(1-\eta^{2})H}}\int^{t}_{t_{0}}\frac{1}{|u_{R}(s)|^{2}}\frac{d}{ds}|u_{R}(s)|ds\nonumber\\
&\leq& \mbox{}\left(A(t_{0})+
\frac{M_{10}}{|u_{R}(t_{0})|^{\beta-1}}\right)\frac{1}{\sqrt{2(1-\eta^{2})H}}\frac{1}{|u_{R}(t_{0})|}.
\end{eqnarray*}
By energy equation and the definition of $t_{0}$, we have
\begin{eqnarray*}
A(t_{0})=\omega(t_{0})|u_{R}(t_{0})||\dot{u}_{R}(t_{0})|\leq
\eta|u_{R}(t_{0})|\sqrt{2(H-V(u_{R}(t_{0})))},
\end{eqnarray*}
which implies that for some $M_{8}$, $M_{9}>0$ independent of $R$
\begin{eqnarray*}
\left|\frac{u_{R}(t)}{|u_{R}(t)|}-\frac{u_{R}(t_{0})}{|u_{R}(t_{0})|}\right|\leq
M_{8}\eta+\frac{M_{9}}{|u_{R}(t_{0})|^{\beta}},
\end{eqnarray*}
which proves this lemma.

Since we have Theorems 6.1-6.3, similar to \cite{3}, we have the
following theorem.

\vspace{0.3cm}{\bf Lemma 6.4(See\cite{3})}\ {\em For any
$\varepsilon>0$, there exists $M_{11}>0$ such that for $R>M_{11}$}
\begin{eqnarray*}
u_{R}\left(\left[t^{*},\frac{T_{R}}{2}\right]\right)\bigcap\{|x|\geq
M_{11}\}\subset\left\{y\in
R^{N}:\left|\frac{y}{|y|}-e\right|<\varepsilon\right\},
\end{eqnarray*}
where $e$ is the given direction defined in $E_{R}$ and $t^{*}$ is
defined as (\ref{31}).

Let $\overline{t}\geq t^{*}$ such that
$|u_{R}(\overline{t})|=\overline{L}_{\eta}$. Then we can get for any
$\varepsilon>0$
\begin{eqnarray}
\left|\frac{u_{R}(t)}{|u_{R}(t)|}-e\right|<\varepsilon
\end{eqnarray}
for all $t\geq\bar{t}$, which implies that
\begin{eqnarray*}
\frac{u_{\infty}(t)}{|u_{\infty}(t)|}\rightarrow e\ \ \ \mbox{as}\ \
\ t\rightarrow+\infty
\end{eqnarray*}
and
\begin{eqnarray*}
\frac{u_{\infty}(t)}{|u_{\infty}(t)|}\rightarrow e\ \ \ \mbox{as}\ \
\ t\rightarrow-\infty.
\end{eqnarray*}

From the above discussion, we have proved there is at least one
hyperbolic solution for $(1)-(2)$ with $H>0$ which has the given
asymptotic direction at infinity. We finish the proof. $\Box$

\section*{Acknowledgements} The authors sincerely thank the referees for his/her many valuable
comments and remarks which make the paper more simple and clear.

\end{document}